\documentclass[12pt,a4paper]{amsart}
\usepackage[cp1251]{inputenc}
\usepackage[english,russian]{babel}
\usepackage{amsfonts}
\usepackage{amssymb}
\usepackage{latexsym}
\small\large
\newtheorem{lemma}{Лемма}
\newtheorem{remark}{Замечание}

\newtheorem{proposition}{Предложение}

\def\m{\smallsetminus}

\title {О группах с сильно вложенной унитарной подгруппой }
\author{ А.И. Созутов}
\begin{document}

\maketitle
\sloppy

 \begin{center}
Исследование выполнено за счет гранта Российского научного фонда (проект № 19-71-10017)
\end{center}
\begin{center}
Созутов Анатолий Ильич, Сибирский федеральный университет, пр. Свободный, 79, Красноярск 660041, Россия, E-mail: sozutov\_ai@mail.ru
\end{center}

\begin{abstract}
The proper subgroup $B$ of the group $G$ is called {\it strongly embedded},
if $2\in\pi(B)$ and $2\notin\pi(B \cap B^g)$ for any element
$g \in G \setminus B $ and, therefore, $ N_G(X) \leq B$ for any 2-subgroup
$ X \leq B $. An element $a$ of a group $G$ is called {\it finite} if
for all $ g\in G $ the subgroups $ \langle a, a^g \rangle $ are finite.

In the paper, it is proved that the group with finite element of order $4$ and strongly embedded
subgroup isomorphic to the Borel subgroup of $U_3(Q)$ over a locally finite
field $Q$ of characteristic $2$ is locally finite and isomorphic to the group $U_3(Q)$.

A strongly embedded subgroup of a unitary type, subgroups of Borel, Cartan, involution, finite element.
\end{abstract}

\bigskip

\section*{Введение}

Собственная подгруппа $B$ группы $G$ называется {\it сильно вложенной},
если $2\in \pi (B)$ и
$2\notin\pi (B\cap B^g)$ для любого элемента
$g\in G\setminus B$ и, значит,  $N_G(X)\leq B$ для любой 2-подгруппы
$X\leq B$ \cite[определение 4.20]{Gor}. Как пишет Д. Горенстейн
\cite[стр. 26-27, 196-202]{Gor}, понятие сильно вложенной подгруппы
составляет один из наиболее важных инструментов теории конечных простых групп.
Согласно известным  результатам М. Судзуки и Х. Бендера \cite{SM,BH}
конечная простая  группа с сильно вложенной подгруппой изоморфна одной из
групп $L_2(2^n)$, $Sz(2^n)$, $U_3(2^n)$ \cite[теорема 4.24]{Gor}.

Элемент  $a$ группы $G$ называется {\it конечным}, если
для всех $g\in G$ подгруппы $\langle a,a^g\rangle$ конечны.
Так, например, в периодической группе каждая инволюция является конечным
элементом. В \cite{S,SS1} было доказано, что группа с конечной инволюцией и
сильно вложенной подгруппой, изоморфной подгруппе
Бореля группы $L_2(Q)$ или $Sz(Q)$ над локально конечным полем $Q$ характеристики $2$,
локально конечна и изоморфна соответственно группе $L_2(Q)$ или $Sz(Q)$.
Эти результаты  были востребованы и развиты в исследованиях групп с условиями насыщенности
\cite{SS2} -- \cite{LS2}, в том числе  были получены и характеризации
унитарных групп $U_3(Q)$ из списка Судзуки-Бендера.
В настоящей  работе
получена характеризация групп $U_3(Q)$  без условия насыщенности:

\bigskip

\noindent{\bf Теорема. }{\it
Группа с конечным элементом порядка 4 и
сильно вложенной подгруппой, изоморфной подгруппе
Бореля группы $U_3(Q)$ над локально конечным полем $Q$ характеристики $2$,
локально конечна и изоморфна группе $U_3(Q)$.
}

\smallskip

\section{Определения и используемые результаты}

В доказательстве  используется  теорема
О.Н. Кегеля и Б.А. Верфрица \cite{KegW}:

\begin{proposition} \label{p1} 
Если локально конечная группа
обладает локальной системой, состоящей из возрастающей цепочки
конечных простых групп лиева типа ранга 1, то она изоморфна
некоторой простой группе лиева типа ранга 1.
\end{proposition}

Согласно этой теореме  исследуемая группа $G$
будет представлена как объединение $G=\bigcup_{i=2}^{\infty} G_i$
найденной в лемме \ref{l11} возрастающей цепи конечных подгрупп
\begin{equation}
G_2<G_3<...\, ,  \ \ \mbox{где } G_i=U_3(Q_i),\  Q_i\mbox{ --- конечное поле
характеристики 2.}
\label{f01}
\end{equation}

В бесконечной локально конечной группе $G=U_3(Q)$ множество цепей вида
(\ref{f01}) имеет мощность континуума, однако строение $G$ и ее подгруппы Бореля
$B=U\lambda H$ однозначно определяется строением поля $P$  c автоморфизмом $\sigma$ порядка 2 и
его подполя  $Q$ неподвижных точек относительно $\sigma$.
В работе используются обозначения
$^2A_2(2^{2n})=PSU_3(2^n)=U_3(2^n)$ из \cite[таблица
6.4.1]{Kon}, \cite[стр. 166-168]{Gor}, и в формулировках утверждений
используется поле $Q$, а основное поле $P$ "остается за кадром". Надеюсь, что
это не вызовет недоразумений, тем более что в используемых в работе
результатах из \cite{Gor} применяются такие же обозначения.

Локально конечное поле $Q$ счетно и
является объединением $Q=\cup_{i=2}^{\infty} Q_i$ конечных подполей $Q_i$,
составляющих цепь $Q_2\subset Q_3\subset...$ .
Каждой такой цепи соответствует возрастающая цепь $B_2<B_3<...\,$ подгрупп
Бореля $B_i=U_i\leftthreetimes H_i$ групп
$G_i=U_3(Q_i)$, объединение $B=\cup_{i=2}^{\infty} B_i$  которой
совпадает c заданной условиями теоремы сильно вложенной подгруппой.
Силовские 2-подгруппы $U_i$ групп $B_i$ также составляют цепь
$U_2<U_3<...\,$, $U=\cup_{i=2}^{\infty} U_i$; аналогично
задана цепь подгрупп Картана $H_2<H_3<...\,$, $H=\cup_{i=2}^{\infty} H_i$.
По п. 2 предложения \ref{p2} $H=H_0\times H_1$ --- локально циклическая
группа, $H_0$ изоморфна мультипликативной группе поля $Q$ и является
объединением $H_0=\cup_{i=2}^{\infty} H_{0i}$ цепи подгрупп $H_{0i}\leq H_i$, изоморфных мультипликативным группам полей $Q_i$
\begin{equation}
H_{02}<H_{03}<...\,,\  \mbox{ где } H_{0i}=H_0\cap H_i.
\label{f1}
\end{equation}
Цепь (\ref{f1}) однозначно определяет цепь $Q_2\subset Q_3\subset ...$
подполей $Q_i$. Отметим этот момент, поскольку в доказательстве
теоремы  движение происходит в обратном  изложенному выше направлении,
от цепи (\ref{f1}) к цепи (\ref{f01}).

В доказательстве леммы
\ref{l4} используется теорема 4.1 из \cite{S}:

\begin{proposition}
Пусть группа $G$ содержит конечную инволюцию и
элементарную абелеву
$2$-подгруппу $Z$, нормализатор $N_G(Z)$ сильно вложен в $G$ и
является группой Фробениуса с ядром $Z$ и
локально циклическим дополнением $T$. Тогда $G$ локально конечна и
изоморфна группе $L_2(Q)$, где $Q$ ---
локально конечное поле характеристики $2$.
\label{p1}
\end{proposition}

Нам понадобятся свойства подгрупп Бореля групп $U_3(Q)$, \cite[стр. 166-168]{Gor}:

\begin{proposition}\label{p2} Пусть $Q$ ---  (локально) конечное поле
характеристики $2$, \linebreak  $B=U\leftthreetimes  H$ --- подгруппа Бореля
группы $G=U_3(Q)$, $H$ --- ее подгруппа Картана, $U$ --- силовская
$2$-подгруппа групп $B$ и $G$,  и $Z=Z(U)$. Тогда
\begin{enumerate}
\item $B=U\leftthreetimes  H,\ H$ --- (локально) циклическая группа,
$U$ --- силовская 2-подгруппа в $G$  ступени нильпотентности $2$,
периода $4$ и $Z=Z(U)=U'=\Phi
(U)=\Omega_1(U)$;
\item $H=H_0\times H_1$, где $H_1=C_H(Z)$,
$H_0$ и $H_1$  --- (локально) циклические подгруппы,   пересечение $\pi (H_0)\cap\pi
(H_1)$ пусто;
\item Если $|Q|$ --- конечное поле из $q=2^n$ элементов, то
$|G|=q^3(q^2-1)(q^3+1)/d$,
$|U|=q^3$, $|Z|=q$, $|H|=(q^2-1)/d$,
$|H_1|=(q+1)/d$, $|H_0|=q-1$, где $d=(3,q+1)$;
\item $U\leftthreetimes  H_0$ --- группа Фробениуса с неинвариантным множителем
$H_0$, действующем транзитивно на множестве $Z^\#$ всех инволюций группы $U$;
\item $B/Z$ --- группа Фробениуса с ядром $U/Z$
и дополнением $HZ/Z$, и фактор-группа $HZ/Z$
либо транзитивна на множестве неединичных элементов фактор-группы $U/Z$,
когда в $H_0$ есть элемент порядка $3$, либо имеет ровно
3 орбиты.
\end{enumerate}
\end{proposition}

Элементарные свойства произвольной группы  $G$
с конечной инволюцией и сильно вложенной подгруппой $B$
можно найти в \cite[леммы 1.7, 2.1, 2.2]{S} и в
\cite{SS1,SSS}:

\begin{lemma}\label{l1} Пусть $G$ --- группа  c конечной инволюцией и
сильно вложенной подгруппой $B$, $J$ --- множество инволюций в $G$,
$i\in J\cap B$, $k\in J\setminus B$ и $K=B\cap B^k$. Тогда
\begin{enumerate}
\item Все инволюции в группе $G$ сопряжены и конечны, $J\cap B=i^B$,
порядок произведения $ik$ нечетен и множество $J$ инволюций в $G$
равномерно распределено по смежным классам
$G/B$: $|J\cap Bg|=|J\cap B|$ для любого элемента $g\in G$;
\item Если некоторая
силовская 2-подгруппа $S$ группы $G$ имеет нетривиальное
пересечение с  $B$, то она содержится в $B$;
\item Если  $L<G$, $L\nleq B$ и  $2\in\pi (L\cap B)$, то
подгруппа $L\cap B$ сильно вложена в $L$;
\item Если $b$ --- неединичный элемент нечетного порядка из $G$ и $b^j=b^{-1}$
для некоторой инволюции $j\in G$, то в централизаторе $C_G(b)$ нет инволюций;
\item  В каждом смежном классе
$C_B(i)\cdot b$,
где $b\in B$, существует единственный инвертируемый инволюцией $k$
элемент и его порядок конечен и нечетен;
\item $B=K\cdot C_B(i)=T\cdot C_B(i)$ и $i^B=i^T$,
где $T$ --- подгруппа в $K=B\cap B^k$, порожденная инвертируемыми
инволюцией $k$ элементами из $K$.
\end{enumerate}
\end{lemma}

\section{\bf Леммы о группе с конечной инволюцией}

В данном параграфе $G$ --- группа с конечной инволюцией
и сильно вложенной подгруппой $B$, изоморфной подгруппе
Бореля группы $U_3(Q)$ над локально конечным полем $Q$ характеристики $2$.
Здесь и далее в тексте $B=U\leftthreetimes  H$,   $U$, $H$, $H_0$,  $H_1$ и $Z$ --- подгруппы
группы $G$, определенные
в предложении \ref{p2}.

\begin{lemma}\label{l2}
Пусть $k$ ---  инволюция из $G\m B$ и
$T$ --- подгруппа, порожденная всеми инвертируемыми инволюцией
$k$ элементами из $B\cap B^k$. Тогда $T=H_0^u$ для подходящего
элемента $u\in U$, в частности,  $H_0$ инвертируется некоторой
инволюцией $v\in G\m B$.
\end{lemma}
\begin{proof}
По условиям в подгруппе  $K=B\cap B^k$ нет инволюций,
 согласно п. 1 предложения \ref{p2} $K\simeq KU/U$ ---
локально циклическая группа и множество инвертируемых инволюцией
$k$ элементов из $K$ является подгруппой $T=\{ t\in K\,|\, t^k=t^{-1}\}$.
По п. 4 леммы \ref{l1} $C_U(t)=1$
для любого $t\in T^\#$. Поэтому $T$ действует сопряжением свободно
на $U$.  Группа $B$ локально конечна,  для ее подгрупп верна теорема
Фробениуса и, значит,   подгруппы
$Z\leftthreetimes T$, $U\leftthreetimes T$  являются
группами  Фробениуса с дополнением $T$.
Ввиду п. 1 предложения \ref{p2} в локально конечной группе
$B$ все конечные $p$-подгруппы для нечетных $p$ циклические, и
по теореме Силова все подгруппы из $B$ одного и того же нечетного простого
порядка $p$ сопряжены. По  п. 2 предложения \ref{p2}
$H_1\leq C_B(Z)$, следовательно $\pi (T)\subseteq \pi (H_0)$
и $T$ содержится в группе Фробениуса $U\leftthreetimes H_0$
(п. 4 предложения \ref{p2}). В силу свойств (локально) конечных групп
Фробениуса \cite[предложение 1.14]{PSS} $T$ содержится в одном из дополнений
группы $U\leftthreetimes H_0$, и с точностью до сопряженности в $U\leftthreetimes H_0$
можно считать,  что  $T\leq H_0$.
В силу пп. 5, 6 леммы \ref{l1} $T$ действует транзитивно на $Z^\#$
и из п. 4 предложения \ref{p2} следует равенство $T=H_0$.
Как показано выше, инволюция $v=k$ инвертирует $T$.
Лемма доказана.
\end{proof}

\smallskip

Зафиксируем обозначение инволюции $v$ из леммы \ref{l2} до конца работы.

\begin{lemma}\label{l3}
Подгруппа $H$ совпадает с централизатором $C_G(t)$ каждого
неединичного элемента
$t\in H_0$. Справедливы равенства $C_B(v)=H_1$,
$B\cap B^v=H$  и $N_G(H)=H\leftthreetimes\langle v\rangle$.
\end{lemma}

\begin{proof}
Пусть $1\ne t\in H_0$ и предположим,
что $C_G(t)\ne H$. В силу пп. 1, 4 предложения \ref{p2} $C_B(t)=H$,
следовательно в $C_G(t)\setminus B$ найдется  элемент $g$.
По п. 1 леммы \ref{l1} $g=bk$, где $b\in B$,
$k$ --- инволюция из $G\setminus B$, и
 $t\in K=B\cap B^k$. По лемме \ref{l2} инвертируемые инволюцией
$k$ элементы из $K$ составляют подгруппу $T=H_0^u$ и ввиду п. 4
предложения \ref{p2} $U\leftthreetimes T=U\leftthreetimes H_0$ --- группа
Фробениуса с  ядром $U$ и дополнениями $T$ и $H_0$.
Поскольку $K$ не содержит инволюций, то $K\cap U=1$, $t\in T\cap H_0$ и
$T=H_0$ как дополнения одной группы Фробениуса.
По лемме \ref{l2} $t^k=t^{-1}$ и, значит, $t^b=t^{-1}$,
что противоречит строению подгруппы $B$ (предложение \ref{p2}).
Следовательно, $C_G(t)=H$ для любого неединичного элемента $t\in H_0$.
По лемме \ref{l2} $t^v=t^{-1}$ и, значит, $H^v=H$ и $H\leq B\cap B^v$.
Так как любая подгруппа из $B$ содержащая $H$ собственным образом
содержит инволюции (п. 1 предложения \ref{p2}),
а подгруппа
$K=B\cap B^v$ инволюций не содержит, то $B\cap B^v=H$.
По предложению \ref{p2} $H=H_0\times H_1$ --- локально циклическая группа,
$\pi (H_0)\cap \pi (H_1)$ пусто и $H_1\leq C_B(Z)$.
Отсюда выводим $H_1^v=H_1$ и $H_1=C_H(v)$.

Ввиду пп. 4, 5 предложения \ref{p2} $N_B(H)=H$. Пусть $g\in N_G(H)\setminus
B$. По п. 1 леммы \ref{l1} $g=bk$, где $b\in B$,
$k$ --- инволюция из $G\setminus B$, и $H=B\cap B^g=B\cap B^k$.
Значит, $H^k=H$, $H^b=H$ и $b\in H$.
По лемме \ref{l2} $H_0=\{ h\in H\mid h^k=h^{-1}\}$
и $kv\in C_G(H_0)=H$. Следовательно, $k\in Hv$,  $N_G(H)=H\leftthreetimes \langle v\rangle$
и лемма доказана.
\end{proof}

\begin{lemma}
$C_G(H_1)=H_1\times L$, где  $L=H_0Z\langle v\rangle
Z$ --- подгруппа,  изоморфная $L_2(Q)$.
Для любого элемента $t\in H_1^\#$ справедливо равенство
$N_G(\langle t\rangle)=C_G(H_1)$.
\end{lemma}

\begin{proof}
Рассмотрим фактор-группу $\bar G_0=C_G(H_1)/H_1$ и в ней подгруппу
$\bar B_0=ZH/H_1$. Так как по лемме \ref{l3} $\bar{v}\in \bar{G_0}  $, то $\bar{B_0}\ne\bar{G_0}$.
Пусть $\bar x\in \bar G_0\setminus \bar B_0$ --- произвольный элемент и
$g$ --- один из его прообразов в $G$. Очевидно, что $g\notin B$
и так как пересечение
$ZH\cap (ZH)^g\leq B\cap B^g$ не
содержит инволюций, а $H_1$ --- периодическая группа, то и
$\bar B_0\cap \bar B_0^{\bar x}$ не
содержит инволюций. Следовательно, $\bar B_0$
сильно вложена в $\bar G_0$. Согласно пп. 2, 4 предложения \ref{p2}
$C_G(H_1)\cap B=H_1\times (Z\leftthreetimes H_0)$ и, значит,
$\bar B_0\simeq \bar Z\leftthreetimes \bar H_0$ --- группа
Фробениуса с элементарным абелевым ядром $ZH_1/H_1\simeq Z$
и локально циклическим дополнением $H/H_1\simeq H_0$.
В силу п. 1 леммы \ref{l2} каждая инволюция в группе $\bar G_0$
конечна.
  По предложению \ref{p1}
$\bar G_0\simeq L_2(R)$ для подходящего локально конечного поля $R$
характеристики 2. Так как мультипликативные группы локально конечных
полей $R$ и $Q$ изоморфны $H_0$, то поля $R$ и $Q$ изоморфны
и $\bar G_0\simeq L_2(Q)$. По лемме \ref{l3} $\bar v\in \bar G_0$ и
$\bar G_0=\bar B_0\langle\overline{v}\rangle
\bar B_0$.

По теореме Шмидта \cite{KarM} группа $C_G(H_1)$
локально конечна. Поле $Q$ есть объединение $Q=\cup Q_i$
возрастающей цепочки $ Q_2\subset Q_3\subset\ldots$ конечных полей характеристики 2,
которой в $G_0$
соответствует возрастающая цепочка
\begin{equation}
\overline   K_2<\overline K_3<\dots
\label{f5}
\end{equation}
конечных подгрупп $\overline  K_i\simeq L_2(Q_i)=L_2(2^{n_i})$,
причем $\cup \overline K_i=\bar G_0$.
Так как мультипликатор Шура в $L_2(2^n)$
равен $Z_2$ (\cite{Gor}, табл. 4.1),
а подгруппа $H_1$ не содержит инволюций,
то полный прообраз $K_i$ подгруппы $\overline K_i$
разлагается в прямое произведение $K_i=H_1\times L_i$,
где $L_i\simeq L_2(Q_i)=L_2(2^{n_i})$. Отсюда выводим, что
подгруппы $L_i$ составляют цепочку
\begin{equation}
 L_2< L_3<\dots ,
\label{f6}
\end{equation}
объединение $L$ которой по предложению \ref{p1} изоморфно $L_2(Q)$.
Таким образом,
$C_G(H_1)=H_1\times L$, очевидно $L=ZH_0\langle v\rangle Z$  и первое утверждение леммы доказано.

Пусть $t\in H_1^\#$, $g\in N_G(\langle t\rangle)$ и $g\notin C_G(H_1)$.
По лемме \ref{l1} $g=bk$, где $b\in B$, $k\in J\setminus B$,
а в силу лемм \ref{l2}, \ref{l3}
$\langle t^k\rangle=\langle t\rangle$ и $k\in C_G(H_1)$.
Значит, $b\notin N_B(\langle t\rangle )=H_1\times ZH_0$
и ввиду предложения \ref{p2}
 можно считать,
что $b=u$ --- элемент порядка 4. Но тогда
$\langle t^{uk}\rangle=\langle t\rangle$ и
$\langle t^u\rangle=\langle t\rangle$, $u^t=u$, что противоречит
действию $H_1$ на фактор-группе $U/Z$ (п. 5 предложения \ref{p2}).
Следовательно, $N_G(\langle t\rangle)=C_G(H_1)$
и лемма доказана.
\end{proof}

\begin{lemma}\label{l5}
Подгруппа $UH_0$ действует (сопряжениями) транзитивно
 на множестве инволюций $J\setminus B$, группа $G$ действует дважды
транзитивно на множестве $\Omega =B^G$ подгрупп, сопряженных с $B$,
и подгруппа $U$ действует на $\Omega\setminus\{ B\}$ регулярно.
В частности, $G=B\langle v\rangle B=B\cup BvU$, $U$ --- силовская $2$-подгруппа группы
$G$, силовские $2$-подгруппы
в $G$ сопряжены и попарно взаимно просты.
\end{lemma}

\begin{proof}
По предложению \ref{p2} $H=C_B(H_0)=N_B(H_0)$ и согласно лемме \ref{l3}
и $N=N_G(H_0)=H\leftthreetimes\langle v\rangle$.
По предложению \ref{p2} $UH_0$ --- группа Фробениуса с ядром $U$ и дополнением
$H_0$ и, значит, каждое дополнение $K$ подгруппы $U$ в $B$
имеет вид $K=H^u$, где $u\in U$.
Пусть $k$ --- произвольная инволюция из $J\setminus B$ и
$T$ --- подгруппа, порожденная всеми инвертируемыми инволюцией
$k$ элементами из $B\cap B^k$. По лемме \ref{l2}  $T=H_0^{u^{-1}}$ для подходящего
элемента $u\in U$,  ввиду леммы \ref{l3} $B\cap B^k=H^{u^{-1}}$
и $k^u\in N$.
Из леммы \ref{l3} следует $J\cap N=v^{H_0}$, значит,
$k^{uh}=v$ для подходящего $h\in H_0$. Поскольку $C_G(v)\cap UH_0=1$, то
$UH_0$ действует сопряжениями регулярно
на множестве $J\setminus B$, то есть транзитивно и без неподвижных точек.

Если $V\in Syl_2(G)$ и $V\cap U\ne 1$, то
по п. 2  леммы \ref{l1}  $V\leq B$ и согласно предложению \ref{p2} $V=U$.
Значит, $U\cap V=1$ при $V\ne U$.
Положим $V=U^v$. Из транзитивности действия группы $UH_0$
 на множестве $J\setminus B$ следует
$J\setminus U\subseteq \cup_{x\in UH_0} V^x$ и
$Syl_2(G)\setminus\{ U\}=V^{UH_0}$. Итак, $Syl_2(G)=U^G$
и $G$ действует дважды транзитивно на
множестве $\Omega=B^G$. Так как $V^{H_0}=\{ V\}$ и $N_U(V)=1$,
то $Syl_2(G)\setminus\{ U\}=V^U$ и $U$ действует  сопряжениями на
множестве $\Omega\setminus\{ B\}$ регулярно.
Из доказанного очевидно следует, что $G=B\langle v\rangle B=B\cup BvU$, $U$ --- силовская $2$-подгруппа группы
$G$, силовские $2$-подгруппы
в $G$ сопряжены и попарно взаимно просты.
 Лемма доказана.
\end{proof}

\begin{lemma}\label{l6}
Каждый элемент $g\in G\setminus B$ имеет
единственное каноническое
представление вида $g=hu_1vu_2$, где
$u_1,u_2\in U,\ h\in H$.
Каждая инволюция $t$ из $G\setminus B$  канонически представима в
видe $t=u_t^{-1}h_tvu_t=h_tu_t^{-h_t}vu_t$, где $u_t\in U$, $h_t\in H_0$.
\end{lemma}

\begin{proof}
По лемме \ref{l5} $G=B\cup BvU$ и каждый элемент $g\in G\setminus B$
представим в виде $g=hu_1vu_2$, где
$u_1,u_2\in U,\ h\in H$. Согласно \cite{Gor}[стр. 160] такое
представление называем {\it каноническим}.
Оно единственно, поскольку
из $g=hu_1vu_2=h_1u_3vu_4$ вначале следует
$u_3^{-1}h_1^{-1}hu_1=vu_4u_2^{-1}v$, затем
$u_4=u_2$ ввиду сильной вложенности $B$ в $G$, и, наконец,
$h=h_1$ и  $u_3=u_1$ в силу равенства $B=U\leftthreetimes H$
(предложение \ref{p2}).

По леммам \ref{l5}, \ref{l3} $J\setminus B=v^{H_0U}=(H_0v)^U$ и  инволюция
$t$ из $G\setminus B$  канонически представима в
видe $t=u_t^{-1}h_tvu_t=h_tu_t^{-h_t}vu_t$, где $u_t\in U$, $h_t\in H_0$.
Лемма доказана.
\end{proof}

\bigskip

Воспользуемся обозначениями $Q_k$,  $U_k$,   $H_k$,  $B_k$, здесь
$k\geq 2$,
введенными  в \S 1;
обозначим также $Z_k=Z(U_k)$, $H_{1k}=H_1\cap H_k$ и $H_{0k}=H_0\cap H_k$,
где $H_1$, $H_0$ --- подгруппы из $H$ (предложение \ref{p2}),
так что $H_k=H_{1k}\times H_{0k}$.
 Введем конечные
подмножества
\begin{equation}
M_k=B_k\langle v\rangle B_k=B_k\cup B_kvU_k.
\label{f8}
\end{equation}

Так как $H$ --- локально циклическая группа
(п. 1 предложения \ref{p2}) и $v\in N_G(H)$ (лемма \ref{l3}),
то все подгруппы $H_k$
и подгруппы $H_{1k}$,  $H_{0k}$ допустимы относительно 
инволюции $v$ (леммa \ref{l2}).
Поэтому в (\ref{f8}) выполняется равенство $B_kvB_k=B_kvU_k$.

В группе $Z=Z(U)$ есть однозначно
определенная инволюция $u_0$ \cite{Gor}[стр. 163], для которой выполняется
{\it структурное уравнение Судзуки групп $L$ и $G$}:
\begin{equation}
vu_0v=u_0vu_0, \ \mbox{или}\ (vu_0)^3=1.
\label{f9}
\end{equation}

\begin{remark}
Если $u_1\in Z$,  $u_1\ne u_0$ и $vu_1v=u_1vu_1$, то подгруппа
$\langle u_0,v,u_1\rangle$  изоморфна симметрической группе $S_4$ \cite[\S
 6.2]{CoxM}. Однако в $S_4$ нет сильно вложенных подгрупп,
противоречие п. 3 леммы \ref{l1}.
Следовательно инволюция $u_0$ из (\ref{f9}) единственна в $Z$.
\label{z2}
\end{remark}

Очевидно можно считать, что $u_0\in U_k$ для всех $k$.
Вместе с формулами
умножения в подгруппах
$Z\leftthreetimes H_{0}$ и $H_{0}\leftthreetimes\langle v\rangle$
уравнениe (\ref{f9})
однозначно определяет умножение в группе
$L=ZH_{0}\langle v\rangle Z$. Так как
$H_0v=vH_0$ ($hv=vh^{-1}$), то для любых $z_1,z_2,z_3,z_4\in Z$,
$h_1,h_2\in H_0$ имеем  $(z_1h_1vz_2)\cdot (z_3h_2vz_4)=
(z_1h_1)\cdot (vz_2z_3v)\cdot (h_2^{-1}z_4)$ и с помощью (\ref{f9})
находим представление элемента $vtv$ в виде $z_5h_3vz_6$, здесь $t=z_2z_3\in Z$.
Так как инволюция $t\in Z$ представима в виде
$t=h_t^{-1}u_0h_t$ при подходящем и единственном $h_t\in H_0$
(п. 4  предложения \ref{p2}), то применяя (\ref{f9}), получаем
\begin{equation}
vtv=vh_t^{-1}u_0h_tv=h_tvu_0vh_t^{-1}=
h_tu_0vu_0h_t^{-1}=u_0^{{h_t}^{-1}}h_t^2vu_0^{h_t^{-1}},
\label{f13}
\end{equation}
где $z_5=u_0^{{h_t}^{-1}}$, $h_3=h_t^2$, $z_6=u_0^{h_t^{-1}}$.
Тем самым
умножение в $L$ определено однозначно
и  установлен изоморфизм
$L\simeq L_2(Q)$.  Те же вычисления верны и когда
$L=L_k$ --- конечное множество, из них следует, что $L_k$ --- подгруппа и
изоморфизм $L_k\simeq L_2(Q_k)$. Итак, верна следующая

\begin{lemma}\label{l7}
Множество $M_k$ содержит подгруппы $L_k=Z_kH_{0k}\langle v\rangle Z_k$
и $H_{1k}\times L_k$.
\end{lemma}

\section{\bf Доказательство теоремы}

\begin{lemma}\label{l8}
Пусть $u$ --- элемент порядка $4$ из $U$ и подгруппа
$M=\langle u,v\rangle$ конечна. Тогда либо $M=A\leftthreetimes \langle
u\rangle$ группа Фробениуса с дополнением $\langle
u\rangle$ и абелевым ядром $A$, либо $M$ изоморфна
группе $U_3(Q_m)$ над некоторым конечным подполем $Q_m\subset Q$.

\end{lemma}

\begin{proof} По лемме \ref{l1}  подгруппа $B_m=B\cap M$ сильно вложена в
$M$ и  $U_m=B_m\cap U$ --- нормальная в $B_m$ силовская 2-подгруппа группы $M$.
В силу леммы \ref{l5} силовские $2$-подгруппы
в $M$ попарно взаимно просты, в частности $O_2(M)=1$.

\smallskip

Если инволюция $u^2$ в $U_m$ одна, то  $M=A\cdot C_M(u^2)$, где $A=O(M)$,
при этом $v\in A\leftthreetimes \langle u^2\rangle$ и, значит,
$M=A\leftthreetimes \langle u\rangle$. Следовательно,
$U_m=\langle u\rangle$, ввиду предложения \ref{p2} $A\cap B=1$, $B_m=U_m=\langle
u\rangle$, $C_A(u^2)=1$ и
$M=A\leftthreetimes \langle
u\rangle$ --- группа Фробениуса с дополнением $\langle
u\rangle$ и абелевым ядром $A$, $a^v=a^{u^2}=a^{-1}$ для любого $a\in A$.

\smallskip

Пусть $Z_m=U_m\cap Z$ не циклическая группа. Тогда ввиду предложения \ref{p2}
и леммы \ref{l1}
$B_m$ содержит
подгруппу $Z_m\leftthreetimes H_{0m}$, где $H_{0m}\leq H_{0}$
и $H_{0m}$ действует сопряжением регулярно на $Z^\#_m$ в силу пп. 2, 3
предложения \ref{p2}.
Согласно лемме
\ref{l7} $M$ содержит подгруппу $L_m=Z_mH_{0m}\langle v\rangle Z_m$,
изоморфную $L_2(Q_m)$, где $|Q_m|=|H_{0m}|+1=q_m$ (подполе $Q_m$ в $Q$
совпадает очевидно с множеством решений в $Q$ уравнения $x^{q_m}-x=0$).
Поскольку $U_m$ нормальна в $B_m$ и содержит элемент порядка 4,
то в силу строения подгруппы $B$ (предложение \ref{p2}) $O(B_m)=1$.
Не циклическая группа $Z_m$ не может действовать на конечной группе
без неподвижных точек, поэтому из  $O(B_m)=1$ следует $O(M)=1$.
Итак, $O_2(M)=O(M)=1$ и по лемме \ref{l1} все инволюции в $M$ сопряжены.
Следовательно минимальная нормальная подгруппа $K$ в группе $M$ порождена
множеством $J\cap M$ и является неабелевой простой группой. Так как
$v,u^2\in K$, то $[M:K]\leq 2$. Однако, с другой стороны, $Z_m\leq K$,
$|U_m/Z_m|\geq |H_{0m}|>2$ и ввиду п. 4 предложения \ref{p2} фактор-группа
$B_m/Z_m$ является группой Фробениуса с дополнением нечетного порядка
и в ней нет подгрупп индекса 2. Следовательно $B_m\leq K$ и $K=M$.
Поскольку  $M$ содержит подгруппу
изоморфную $L_2(Q_m)$ (лемма \ref{l7}), то $|М|$ делится на 3. Учитывая
указанные свойства заключаем, что согласно \cite[теорема 4.24]{Gor} $M$
изоморфна $U_3(Q_m)$.
Лемма доказана.
\end{proof}

\begin{lemma}\label{l9}
Множество подгрупп $M_z=\langle uz,v\rangle$,
где $u$ --- фиксированный элемент порядка $4$ из $U$, а $z$ пробегает $Z$,
содержит не более одной конечной непростой группы.
\end{lemma}

\begin{proof} Допустим, что $M=\langle u,v\rangle$  и
$M_z=\langle uz,v\rangle$  --- две различные конечные не простые группы.
По лемме \ref{l8} $M=A\leftthreetimes \langle
u\rangle$ и $M_z=A_z\leftthreetimes \langle uz\rangle$ есть группы
Фробениуса с абелевыми ядрами $A$ и $A_z$ соответственно и различными
дополнениями $\langle u\rangle$,  $\langle uz\rangle$, так что $1\ne z\ne u^2$.
Обозначим $t=u^2=(uz)^2$, $c=vt$ и $d=c^u$.
Имеем $cd=dc$, $c,d^z\in A_z$, $cd^z=d^zc$ и $d,d^z\in C_G(c)$.
Поскольку $c^v=c^{-1}$, то по лемме \ref{l1}
в подгруппе $C_G(c)$ нет  инволюций и  $D=\langle d,d^z\rangle$ ---
подгруппа без инволюций. Так как $d^t=d^{-1}$, $(d^z)^t=(d^z)^{-1}$, то
четверная подгруппа $T=\langle t\rangle\times \langle z\rangle$
содержится в $N_G(D)$,  $K=\langle D,T\rangle=D\leftthreetimes T$
и инволюции $t$, $z$ в $K$ не сопряжены.
Ввиду предложения \ref{p2} $d\notin B$,  $K\nleq B$,
по лемме \ref{l1} подгруппа $K\cap B$ сильно вложена в $K$
и все инволюции в $K$ сопряжены.
Полученное противоречие доказывает лемму.
\end{proof}

\begin{lemma}\label{l10}
Каждое множество $M_k$ из (\ref{f8}) содержится в некоторой конечной
подгруппе $V$ из $G$, изоморфной $U_3(R)$, где $R$ конечное подполе
из $Q$ и $Q_k\subseteq R$.
\end{lemma}

\begin{proof} По условиям теоремы $U$ содержит конечный в $G$ элемент
$u$ порядка 4. В силу п. 4 предложения \ref{p2}  с точностью до сопряженности
в $B$ можно считать, что $u\in U_k$ для всех $k$, и
ввиду транзитивности $H_0$ на $Z^\#$, что $u^2=u_0$.
Исключая тривиальные случаи малых порядков, считаем $|U_2|>2^{4}$.
Пусть $H_{0k}=\langle h_k\rangle$ и $x=u^{h_k}$.
При сделанных предположениях для группы $U_k$ непосредственно проверяется,
что $|x^{U_k}|>4$  и, значит, существует
$z\in Z^\#_k$, $u_0\ne z\ne u_0^{h_k}$, для которой $zx\in x^{U_k}$.
По условию теоремы подгруппы $\langle x,v\rangle$ и $\langle xz,v\rangle$
конечны, по лемме \ref{l9} хотя бы одна из них, обозначим ее $V$, проста
и  изоморфна группе $U_3(R)$, где $R\subset Q$.
Понятно, что $B_m=B\cap V$ --- подгруппа Бореля группы $V$,
$B_m=U_m\leftthreetimes H_m$, $U_m=U\cap V$, $H_m=H\cap V$ и
$H_m=H_{1m}\times H_{0m}$, где $H_{0m}=H_m\cap H_0$,
$H_{1m}=H_m\cap H_1$. При этом ввиду лемм \ref{l2} -- \ref{l4}
$v\in N_{V}(H_m)$ и $V=B_m\cup B_mvU_m$.

Группа $U_m$ содержит инволюцию $u_0^{h_k}=x^2=(xz)^2$ и единственную
 в $Z$ (замечание \ref{z2}) инволюцию $u_0$, удовлетворяющую структурному уравнению
 Судзуки (\ref{f9}). Согласно п. 4 предложения \ref{p2}  инволюции
$u_0$ и $u_0^{h_k}$ сопряжены в $Z_mH_{0m}$ и в $ZH_0$ единственным
элементом $h_k$ из $H_0$, значит, $h_k\in H_{0m}$,
$H_{0k}\leq H_{0m}$ и $Z_k\leq Z_m$ ($Z_k^\#=u_0^{H_{0k}}$).
Множество решений уравнения $x^{q_k}=x$, где $q_k=|Q_k|=|H_{0m}|+1$,
в поле $R$ составляет
подполе, изоморфное $Q_k$. Поэтому $H_{1m}$ содержит единственную
в $H_1$ подгруппу порядка $\frac{q_k+1}{(3,q_k+1)}$ (п. 3 предложения \ref{p2}),  совпадающую
с подгруппой $H_{1k}$. Значит, $H_{1k}\leq H_m$, $H_k\leq H_m$,
и поскольку $H_k$ на $U_k/Z_k$ неприводима, то
$U_k=\langle x^{H_k}\rangle Z_k\leq U_m$.
Итак, $B_k\leq B_m$, $M_k=B_k\langle v\rangle U_k\subseteq V$
и лемма доказана.
\end{proof}

\begin{lemma}\label{l11}
Группа $G$ локально конечна и изоморфна группе $U_3(Q)$.
\end{lemma}

\begin{proof}
Конечные множества
$M_k$, определенные в (\ref{f8}), составляют цепь, а их объединение
$\cup_{k=2}^{\infty} M_k$
в силу   леммы \ref{l6} совпадает с $G$:
\begin{equation}
M_2\subset M_3\subset...\, , \quad G=\bigcup_{k=2}^{\infty} M_k.
\label{f10}
\end{equation}
Согласно лемме \ref{l10}   $M_k\leq V_k<G$,
$V_k\simeq U_3(R_k)$, где $R_k$ конечное поле и $Q_k\subseteq R_k\subset
Q$. Для некоторых из подгрупп $V_k$ могут выполняться равенства
$V_k=M_{k'}$ ($k'\geq k$). Если таких совпадений для членов
$M_{k'}$ цепи (\ref{f10}) бесконечно много (для $k'=k_1,k_2,...$)
то соответствующие подгруппы $V_{k_i}=M_{k_i}\simeq U_3(R_{k_i})$ составляют
цепь $V_{k_1}\leq V_{k_2}\leq ...$ вида \ref{f01} с объединением $\cup_{i=1}^\infty =G$,
удовлетворяющую условиям предложения \ref{p1}. Цепь конечных полей
$R_{k_1}\subseteq R_{k_2}\subseteq ...$ очевидно
является подцепью цепи $Q_2\subset Q_3\subset ...$
и $\cup_{i=1}^\infty R_{k_i}=Q$.
 Согласно предложению
\ref{p1} $G$ локально конечна и изоморфна $U_3(Q)$, и в этом случае лемма доказана.

\smallskip

Если  совпадений $V_k=M_{k'}$ ($k'\geq k$) для  цепи (\ref{f10}) конечное
число, то удалим такие множества $M_{k'}$ из (\ref{f10})
и тогда $V_k\ne M_{k'}$ для всех  $k,k'$.

Укажем цепь $G_{2}<G_{3}<...$
вида (\ref{f01}) конечных подгрупп $G_m\simeq U_3(P_m)$ c соответствующей
цепью $P_2\subset P_3 \subset ...$ конечных подполей поля $Q$.
При $m=2$ полагаем $G_2=V_2\simeq U_3(R_2)$ и $P_2=R_2$.
При $m=3$ находим в цепи (\ref{f10}) множество $M_{k_2}$,
в которой содержится подгруппа $V_2$, и полагаем
$G_3=V_{k_2}$, $P_3=R_{k_2}$ (поскольку $V_{k_2}\simeq U_3(R_{k_2})$ по определению
групп $V_k$).
Аналогично, при $m=n+1$ находим множество $M_{k_{n+1}}$
содержащее подгруппу $G_n$ и полагаем $G_{n+1}=V_{k_{n+1}}$,
$P_{n+1}=R_{k_{n+1}}$,  обеспечивая изоморфизм
$G_{n+1}\simeq U_3(P_{n+1})$ и вложения $G_n<G_{n+1}$,
$P_n\subset P_{n+1}$.

Для  цепи $G_{2}<G_{3}<...$ конечных подгрупп $G_n\simeq U_3(P_n)$
имеем $\cup_{n=1}^\infty G_n=G$, $P_2\subset P_3 \subset ...$,
$\cup_{n=1}^\infty P_n=Q$ и cогласно предложению
\ref{p1} $G$ локально конечна и изоморфна $U_3(Q)$. Лемма и вместе с ней
теорема полностью доказаны.
\end{proof}

\bigskip

Созутов Анатолий Ильич

Сибирский федеральный университет,

пр. Свободный, 79, Красноярск 660041,

E-mail: aisozutov@mail.ru

\end{document}